\documentclass[11pt,reqno]{amsart}
\usepackage{amsmath,amsfonts,amssymb,amsthm,amscd}
\makeindex

\usepackage{amssymb}
\usepackage[latin1]{inputenc}
\usepackage{ epsfig,epsf}
\usepackage{enumerate}
\usepackage{indentfirst}
\usepackage{euscript}
\usepackage{graphicx}
\usepackage{amsmath}
\usepackage{amsfonts}
\usepackage{amssymb}
\makeindex

 \theoremstyle{plain}
 
\newtheorem{proposition}{Proposition}

\newtheorem{remark}{Remark}
\theoremstyle{definition}
 
\newtheorem{definition}{Definition}
\newtheorem{theorem}{Theorem}

\usepackage{amsmath,amsfonts,  amssymb,amsthm,amscd}
\makeindex
\usepackage{amssymb}
\usepackage{graphics}
\usepackage[latin1]{inputenc}
\usepackage{epsf}

\usepackage{ epsfig,epsf}
\usepackage{enumerate}
\usepackage{indentfirst}
\usepackage{euscript}
\usepackage{graphicx}
\usepackage{amsmath}
\usepackage{amsfonts}
\usepackage{amssymb}
\makeindex

 \theoremstyle{plain}

\title[ A Metric Property of Umbilic Points ]{ A  Metric Property
of Umbilic Points } 

 \author[R. Garcia]{Ronaldo Garcia}

 \author[J. Sotomayor]{Jorge Sotomayor}


 \thanks{The first author was partially supported by     FUNAPE/UFG.
Both   authors are fellows of  CNPq.
 This work was done under the project PRONEX/FINEP/MCT - Conv.
76.97.1080.00 - Teoria Qualitativa das Equa\c c\~oes Diferenciais
Ordin\'arias and CNPq - Grant 476886/2001-5.}

\begin{document}
 \maketitle

 \begin{abstract} In the  space $\mathbb U^4$ of 
cubic forms of surfaces, regarded as a $G$-space and endowed with a natural invariant metric,
the ratio of the volumes of those representing umbilic points with negative to those with positive indexes is evaluated in terms of the asymmetry of the metric, defined here. A connection of this ratio with  that   reported by Berry and Hannay (1977) 
in the domain of Statistical Physics, 
is discussed.
  \end{abstract}

 \section{Umbilic Points, Invariant Metrics  and Volume Ratios }

 At an umbilic point $p$ of 
 an oriented $ C^3$ surface $\mathcal S$ embedded
in an oriented Euclidean 3-space ${\mathbb R}^3$ the principal curvatures coincide.  In a neighborhood of such point, $\mathcal S$ can 
be written in a Monge chart  as the graph $z=h(x,y)$ of a function of the form

\begin{equation} \label{eq:1}
 h(x,y) = \frac k2 (x^2+y^2) + \frac 16 (ax ^3+ 3bx ^2y + 3b'xy ^2+ a'y^3) + 
o((x^2+y^2)^{3/2}).
\end{equation}

The frame $(x,y;z)$ is positive and adapted to  $\mathcal S$ at $p$. This means that  the plane orthonormal  frame $(x,y)$ is attached to the tangent plane, positively oriented,  and the $z$-axis is along the unit positive normal to  $\mathcal S$ at $p$. 

 Any other such presentation as the graph  $Z=H(X,Y)$ of a function 

$$ H(X,Y) =\frac K2 (X^2+Y^2) + \frac 16 (AX ^3+ 3BX^2Y + 3B'XY ^2+ A'Y^3) + 
o((X^2+Y^2)^{3/2}$$

 \noindent differs by a rotation 

$$x=\cos{\theta}X-\sin{\theta}Y, \; y=\sin{\theta}X+\cos{\theta}Y, \; z=Z,$$ 

 \noindent linking    the positively oriented frames $(X,Y;Z)$ and $(x,y;z)$,  adapted to the surface at $p$.

The functions are related by $H(X,Y)=h(x,y)$; substitution leads to 

 \begin{equation}\label{eq:2}
\aligned K=&\; k\\
 A =&\;  a\cos^3{\theta}+3b\cos^2{\theta}\sin{\theta}+3b'\sin^2{\theta}\cos{\theta}+a'\sin^3{\theta} \\
B =& -a\sin{\theta}\cos^2{\theta}+b\cos\theta(3\cos^2{\theta} -  2 )+b'\sin\theta (2  -3\sin^2{\theta})+a'\cos{\theta}\sin^2{\theta}\\
 B' =& \;  a\sin^2{\theta}\cos{\theta}+b\sin\theta ( 3\sin^2{\theta}-2)+b'\cos\theta(3\cos^2{\theta} -  2 )+a'\cos^2{\theta}\sin{\theta}\\
A' =&  -a\sin^3{\theta}+3b\sin^2{\theta}\cos{\theta}-3b'\cos^2{\theta}\sin{\theta}+a'\cos^3{\theta} \endaligned
 \end{equation}

 Thus, the group ${\mathbb O}(2)$ of rotations in the plane acts linearly,
to the right, on the four dimensional space  of real cubic forms

$$ \frac 16(ax^3+ 3bx^2y + 3b'xy^2+ a'y^3),$$

\noindent identified with line vectors $u=(a,b,b',a')$ in ${\mathbb R}^4$.

Denote by $\Omega({\theta})$  the  matrix of
the linear transformation in ${\mathbb R}^4$, corresponding to the frame rotation by an angle $\theta$. 
That is $U=u \Omega({\theta})$, with  $U=(A,B,B',A')$ and $u=(a,b,b',a')$.

From equation \ref{eq:2}, get

\begin{equation}
\aligned \Omega(\theta)= \left( \begin{matrix} \cos^3{\theta} &- \sin{\theta}\cos^2{\theta}&\sin^2{\theta}\cos{\theta}& - \sin^3{\theta} \\
 3 \cos^2{\theta}\sin{\theta}  & \cos\theta(3\cos^2{\theta} -  2 ) &- \sin\theta (2  -3\sin^2{\theta})& 3\cos{\theta}\sin^2{\theta}\\
  3 \sin^2{\theta}\cos{\theta} & -\sin\theta ( 3\sin^2{\theta}-2)& \cos\theta(3\cos^2{\theta} -  2 ) & -3\cos^2{\theta}\sin{\theta}\\ 
\sin^3{\theta} &  \sin^2{\theta}\cos{\theta} & \cos^2{\theta}\sin{\theta} &\cos^3{\theta} \end{matrix}\right)\endaligned
 \end{equation}

The space ${\mathbb U}^4$ of umbilic  intrinsic cubic forms on 
surfaces is defined
 as the $G$-Space ${\mathbb R}^4$, endowed  with the above action of the group $G={\mathbb O}(2)$.

The quadratic form 

\begin{equation}\label{eq:in}
T(u)=ab'+a'b-b^2-(b')^2
\end{equation}

\noindent is invariant under $G={\mathbb O}(2)$. That is $T(U)=T(u)$ if $U=u\Omega({\theta})$. Thus it is defined on ${\mathbb U}^4$.  It characterizes the transversal umbilic points, as those with $T(u)\neq 0$. 

It is well known  that the Index $I(u)$ of a transversal umbilic $u$ is $I(u)={\frac12}sign(T(u))$.  See  \cite {bf} and  \cite {gs1, gs2} where the identification of $T\neq 0$ with the transversality to the manifold of umbilic 2-jets is made.

The index of an isolated umbilic counts the number of turns made by  a principal direction at a point of the surface that makes a small circuit around the umbilic, \cite{sp} and \cite{sx}.

 According to \cite{sp} and \cite{st}, the differential equation of principal lines around $p$ in this chart  is defined as a variety in the Projective Bundle. In  the chart $(x,y,[dx:dy])$, the variety is given by the equation:

\begin{equation}\label{eq:lc}
P(x,y;[dx:dy])=Ldy^2 + M dx dy + N dx^2 =0, 
\end{equation}

\noindent where  the functions $L$, $M$ and $N$ are: 
$$\aligned L = h_xh_y h_{yy}- (1+h_y^2)h_{xy} \;\;\;\;\;\;\;\;\; &= - bx -b' y 
 + h.o.t\\
M = (1+h_x^2)h_{yy}- (1+h_y^2)h_{xx} \;\;\; \; &= \;\;(b'-a)x + (a'-b) y+
 h.o.t\\
N = (1+h_x^2)h_{xy} - h_x h_y h_{xx} \;\;\;\;\;\;\;\;\; &=\;\; bx+b' y +  h.o.t 
\endaligned $$

\noindent and therefore  $T(u)=\frac{\partial(N,M)}{\partial(x,y)}|_{(0,0)}$.

\begin{figure}[htbp] 
 \begin{center}
  \hskip .5cm
  \includegraphics[angle=0, width=10cm]{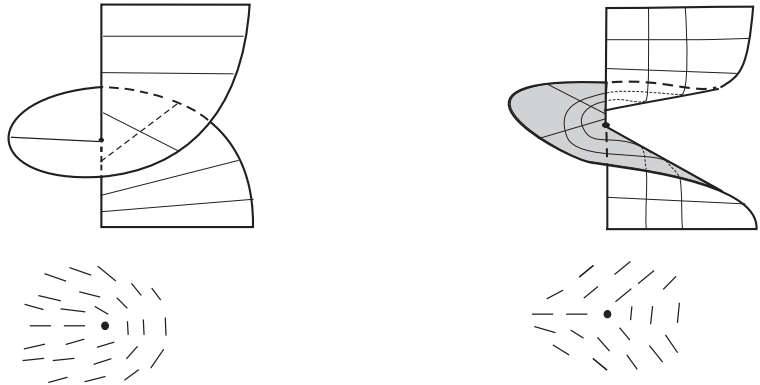}
 \caption{Index of Transversal Umbilic Points: left positive, right negative 
 \label{fig:1}}
 \end{center}
 \end{figure}

\begin{theorem}    \cite{gs1,gs2,bf}.\label{th:1}
Let $p$ be an umbilic point and consider the Monge chart as in equation (\ref{eq:1}).
The transversality condition:  $T \ne 0$ 
holds if and only if 
the surface $P=0$ in equation \ref {eq:lc} is regular along the projective line $x=y=0$.
Then it defines a cylinder  whose projection, removed the projective line, covers a punctured neighborhood of $p$
twice, one for each of the two open cylinders --one for each direction field-- resulting from the removal of the projective line. The covering is orientation preserving or reversing according to $T>0$ or  $T<0$. See Figure \ref {fig:1} for an illustration of one of the cylinders and its projection.
\end{theorem}

A {\it metric} in ${\mathbb U}^4$ is a  positive definite
quadratic form 
such as $q(u)=uQu^{*}$, $Q$ being a positive definite real $4\times 4$ symmetric matrix, invariant under the action of  ${\mathbb O}(2)$.  That is $Q={\Omega({\theta})}Q{\Omega({\theta})}^{*}$, for all ${\theta}$. 

\begin{proposition} \label{pr:1}

Any metric $q(u)=uQu^{*}$, invariant under the action of ${\mathbb O}(2)$, is given by $Q$ of the form:

\begin{equation} \label{eq:mat}
\aligned Q= \left( \begin{matrix} \frac23 \alpha+\frac13 \beta& 0  &\alpha & 0 \\
0 & \beta&0& \alpha\\
  \alpha& 0 & \beta& 0  \\ 0& \alpha &0 &\frac23 \alpha+\frac13 \beta \end{matrix}\right)\endaligned
 \end{equation}

\noindent where $\beta>0, \; \; $
 and $\beta(\frac23 \alpha+\frac13 \beta)-\alpha^2>0$, which gives the positivity of $Q$.
\end{proposition}

\begin{proof}
 Solving the equation $Q={\Omega({\theta})}Q{\Omega({\theta})}^{*}$ first
 for $\theta=\pi/2$ and $\theta=\pi/4$ and finally checking that it holds
 for all $\theta$, the proof follows. 

The criterion for the positivity of a
symmetric matrix,  consisting in that of  all the  principal minors,  is proved in Gantmacher \cite {ga}, Chap. X, pg. 306. 
Notice that for $Q$ in (\ref{eq:mat}) the positivity of the second and third principal minors imply that of the other two.
 \end{proof}

\begin{remark}\label{rem:1}
Another way to obtain the expression of $Q$ in (\ref{eq:mat}) consists in projecting the $10$ dimensional space ${\mathcal M}^{10}$ of $4\times 4$ symmetric matrices  $M$ via the averaging $\mathcal A$ along the orbits of $\Omega(\theta)$:

\begin{equation} \label{eq:av}
{\mathcal A }(M)=\frac{1}{2\pi}\int_0^{2\pi} {\Omega({\theta})}M{\Omega({\theta})}^{*} d\theta .
\end {equation}

Denoting by $m_{ij}$ the entries of the symmetric matrix $M$, integration  in expression (\ref{eq:av}) gives :

\begin{equation} \label{eq:ab}
\begin{aligned}
      \alpha  &=  (6 m_{24} + 3 m_{11} + 3 m_{44} - m_{33} - m_{22} + 6 m_{13})/16\\
       \beta &=   (-6 m_{24} + 9 m_{11} + 9 m_{44} + 5 m_{33} + 5 m_{22} - 6 m_{13})/16
\end{aligned}
\end {equation}
\noindent for the invariant  symmetric matrix ${\mathcal A }(M)$. The other entries of $Q$ in (\ref{eq:mat}) are also  corroborated by integration in (\ref{eq:ab}).

For the identity matrix $I$, ${\mathcal A }(I)$ has  $\alpha=1/4, \; \beta = 7/4$.
\end{remark}

\begin{proposition} \label{pr:2}
 The planes ${\mathbb U}_1 : \; a=3b', a'=3b$ and ${\mathbb U}_2 : \; a= -b', a'=-b$ are invariant under the action of ${\mathbb O}(2)$.
These spaces are mutually orthogonal, relative to $q$.

The quadratic  forms 
$r_1={r_{11}}^2+{r_{12}}^2$ and $r_2={r_{21}}^2+{r_{22}}^2,$
 where 
$$r_{11} = (a+b')/8,\; r_{12} = -(a'+b)/8, \; r_{21} = (a-3b')/24, \;  r_{22} = (a'-3b)/24,$$ 
\noindent are invariant under the action of ${\mathbb O}(2)$. 
Also $r_1$ and $r_2$  vanish respectively on ${\mathbb U}_2$ and  ${\mathbb U}_1$.

Furthermore, the symmetric matrices $R_1$ and $R_2$ which define $r_1$ and $r_2$  generate the lines $\beta=\alpha$ and $\alpha =- \beta /3$, which form the
border of the admissible region $\beta(\frac23 \alpha+\frac13 \beta)-\alpha^2>0$, in Proposition \ref{pr:1}.  
\end{proposition}

\begin{proof}
The invariance of the planes is straightforward.

The plane ${\mathbb U}_1$ is spanned by $u_{1,1}=(3,0,1,0)$ and $u_{1,2}=(0,1,0,3)$; ${\mathbb U}_2$ is spanned by $u_{2,1}=(-1,0,1,0)$ and $u_{2,2}=(0,1,0,-1)$. Scalar multiplication relative to $Q$ of these vectors  ends the proof.

The second and third items follow from a straightforward calculation.
\end{proof}

Although other possibilities exist, in this work the forms $r_1$ and $r_2$ will be used as a reference.

\begin{definition} \label{df:1}
Let $q(u)=uQu^{*}$ be as in Proposition \ref{pr:1}.
Write, $q_1={r_1}/m_1^2$, $q_2={r_2}/m_2^2$, for positive  
 constants $m_1$ and $m_2$, uniquely determined by $q_i=q|_{\mathbb{U}_i} , i=1,2$.
 
The {\it asymmetry} of $q$ is defined  by 
the ratio $\sigma (q)=m_2/m_1$.
\end{definition}

Clearly $\sigma (q)=m_2/m_1$ ranges  over all positive reals. An expression for it in terms of $\alpha, \beta $ has been given in remark \ref{rm:as}.

\begin {theorem} \label{th:in}
Let $T$ be the quadratic form  in equation \ref{eq:in}, giving the index of transversal umbilic points. 

Relative to unit ball $B(1,q)=\{ q(u)\leq 1\}$ of  any 
invariant 
metric  $q$ in ${\mathbb U}^4 $, the ratio of the volume $V_{-}$ of 
the cone $C_{-}$, where $T$ is negative, to that of the volume $V_{+}$ of the   cone $C_{+}$, where $T$ is positive,
 is given by 
$9(\sigma (q))^2$, where $\sigma(q)$ is the asymmetry of $q$, as in  definition \ref{df:1} and remark \ref{rm:as}. 
\end{theorem}

\begin{proof}

Direct  calculation leads to 
$$T=72(-r_2 + r_1 /9 ).$$

 Therefore, in terms of $q_1,\; q_2$, 
$$T=72[ (\frac{m_1} 3)^2 {q_1}-{m_2}^2 {q_2}].$$

 The proof consists in  computing the volume $V_-$   of the solid torus cone

$$C_{-}: q_2 \geq (\frac {m_1}{3m_2})^2 q_1 ,\;\; q\leq 1$$

 and divide it by the volume $V_+$  of 
the solid torus cone 

$$C_{+}: q_2  \leq (\frac {m_1}{3m_2})^2 q_1 , \;\; q\leq 1.$$

Let $v_{i1} ,\; v_{i2}$ be an orthonormal basis of ${\mathbb U}_i ,i=1,2$, relative to $q_i$,  so that they form a positive orthonormal frame, relative to $q$, on ${\mathbb U}^4$. 
 
In  $q$-orthonormal coordinates $(x,y,z,w)$ relative to the frame $v_{i1} ,\; v_{i2}, i=1,\,2$, it follows that 

\begin{equation}\label{eq:diag}
   T=72[(\frac {m_1}{3})^2 (x^2+y^2)-m_2^2 (z^2+w^2)], \;\;\;  q=(x^2+y^2)+(z^2+w^2).
\end{equation}
Let $x=r\cos\theta, \; y=r\sin\theta $ and $z= R\cos\gamma,
\; w=R\sin\gamma$, where 
$0\leq r\leq 1, \;\; 0\leq R\leq 1, \;  0\leq \theta \leq 2\pi, \;\; 0\leq \gamma\leq 2\pi$.

The element of volume $dV$ in the metric $q$ is given by 
$dxdydzdw$

Therefore,  $dV =   r R dr dR d\theta d\gamma$ 
and  
$$V_1= \int_{q\leq 1} dV=  4\pi^2 \int_{0\leq r^2+R^2\leq 1} r R dr dR. $$

Considering  the change of coordinates $r=t \cos\beta, \; R=t\sin\beta$, $0\leq \beta\leq \pi/2$, $0\leq t\leq 1$, 
obtain  
$$ \int_{0\leq r^2+R^2\leq 1} r R dr dR = \int_{0\leq \beta\leq \pi/2, \, 0\leq t\leq 1} t^3 \sin\beta \cos\beta dt d\beta = 1/8$$
 
Therefore, as is well known,  the volume of the unit ball in a four dimensional Euclidean space is given by $V_1=\int_{q\leq 1} dV=  \pi^2 /2$. See Courant-John \cite{cj}, pg. 459.

Take $\tan\beta_0=\frac{m_1}{3m_2}$, the volume of the solid torus cone   $C_+$ is given by , 
$$V_+= 4 \pi^2   \int_0^1\int_0^{\beta_0} t^3 \sin\beta \cos\beta  d\beta dt =\frac{\pi^2}{2} sin^2\beta_0.$$ 

Analogously,   
the volume of the solid torus cone    $C_-$ is equal to 
$$V_- = 4  \pi^2    \int_0^1\int_{\beta_0}^{\frac{\pi}2} t^3 \sin\beta \cos\beta   d\beta dt =\frac{\pi^2}{2}(1- sin^2\beta_0).$$

Direct calculation shows that $\frac{V_-}{V_+}=\frac{(1-\sin^2\beta_0)}{\sin^2\beta_0}=\frac{1}{\tan^2\beta_0}=9(\frac{m_2}{m_1})^2$.
\end{proof}

\begin{remark}\label{rm:as}
In terms  of $Q$, as in equation \ref{eq:mat}, $\sigma (q)$ is calculated as follows:
  \begin {equation}\label{eq:as}
\sigma (q)=m_2/m_1 =\sqrt{\frac{3\alpha+\beta}{3(\beta-\alpha)}}.
\end{equation}
\end{remark}

\begin{proof}
In fact, by the uniqueness of the simultaneous diagonalization of the quadratic forms $q$ and $T$, see \cite{ga}  pg. 314,  equation \ref{eq:diag} implies  that the eigenvalues of matrix $M_T$, of $T$,  relative to $Q$, the matrix of $q$,  are 
$\;72(\frac {m_1}{3})^2$ and $\;-72m_2^2$. 

Separate direct calculation of these  relative  eigenvalues, which are those of the matrix  $M_{T}Q^{-1}$,  gives
$\frac{1}{2(3\alpha +\beta)}$ and $-\frac{3}{2( \beta -\alpha) }$.
                                       
 Equating the ratios of the eigenvalues in both calculations gives   $9(m2/m1)^2$=$3\frac{3\alpha +\beta}{\beta -\alpha}$, which amounts to equation  (\ref{eq:as}).                          

\end{proof}

\section {At the Crossroads of Geometry and Global  Analysis}

 The Geometric local properties of umbilic points, regarded as singularities, have been studied 
focusing  the three following main aspects: 
\begin{enumerate}
\item[i)]
{\em Topological}, related to
the Index  sign of the principal line fields around the umbilic. 
\item[ii)]
{\em Focal}, describing the patterns, Hyperbolic
Elliptic, of normal rays envelopes. This aspect is related to Geometric Optics, Catastrophe Theory and Lagrangian Geometry. See \cite{t,w}.
\item[iii)]
{\em Darbouxian},  which counts the number of principal lines  separatrices  
approaching the umbilic, ($D_1$, $D_2$ or $D_3$) and, more generally, describes locally the foliations by principal lines. 
\end{enumerate}

These aspects and their different types are discriminated
 and  analyzed  in terms of suitable  algebraic
conditions in  the $G$-space ${\mathbb U}^4$.  

 A coherent differential geometric and  topological picture of the set morphology
and inclusion relationships between the different sorts of umbilic types has
been established by Porteous, see 
\cite {p}, and previous reference quoted there. See also Zeeman \cite {z},  for the focal aspect, and Darboux \cite{da}, Sotomayor-Gutierrez \cite{gs1, gsln, gs2} and Bruce-Fidal \cite{bf}, for the Darbouxian types. 

The globalization to the whole surface  of the local analysis  of Darboux, in the context of Structural Stability and Genericity of principal foliations, was carried out in \cite{gs1, gsln, gs2}. 

An additional extension  led  Gutierrez, Garcia, Sotomayor and  others, to expand the study of umbilic points and also principal foliations to surfaces and hypersurfaces in ${\mathbb R}^4$. See \cite{garcia}, \cite{gggt}, \cite{ax}.

Other foliations of interest in  Classical Differential Geometry, such as asymptotic lines  and lines of mean curvature,  defined also by quadratic differential equations similar to(\ref{eq:lc}),  have been studied in \cite{a1}, \cite{ggs}, \cite{m}, \cite{g}, \cite{h}.

 There remain  deep open  problems related to  structure of principal foliations around  isolated umbilic points in smooth surfaces, in the non-transversal case. See Mello-Sotomayor \cite{ms}, Smyth-Xavier \cite{sx} and Ivanov \cite{iv}.

 \section { Umbilic Points  in Random Surfaces}

 On the domain of Statistical Physics, but still connected to 
Geometry and Topology,  Berry and Hannay   \cite{bh}  carried
out a quantitative statistical study of the  proportions in which
the different types of umbilic types are distributed in  random surfaces, such as those modeling  an ocean or a lake. 
An issue here is to study   how the  presence of umbilic points in  a random surface influences the reflection on it of  electromagnetic short waves. 
Although this
 work is more  related to the focal interpretation of umbilic points,  it considers explicitly  also their Darbouxian and Index aspects.

This paper is the outcome of an  initial
attempt to  
provide a
mathematical formulation and a proof, in the tradition of Geometry and Classical Analysis, that could correspond to the conclusions of Berry and Hannay, \cite {bh}, reported in the tradition of  Statistical Physics.

Theorem \ref{th:in} suggests a  disagreement  
 with the report of the calculations in  \cite{bh} which claim  that the statistical ratio is always $ 1$,  disregarding of the statistic anisotropy present  in the evaluation.  The asymmetry of the invariant metric, used to make evaluations in this work,  may be considered as a  geometric   counterpart for  the statistic anisotropy.

Considering  only the local  aspect of surfaces at umbilic points,  this discrepancy may be due to the fact that in the calculations  made in  \cite{bh}, the cubic forms 
are
 regarded as vectors in   ${\mathbb R}^4 $, with a fixed frame,  and not as elements of the  $G$-space  ${\mathbb U}^4 $. The effect of this is that the same umbilic on a surface is counted multiple times, one for each rotated frame.

\vskip .2cm

\author{\noindent Jorge Sotomayor\\Instituto de Matem\'{a}tica e
Estat\'{\i}stica,\\Universidade de S\~{a}o Paulo, \\Rua do Mat\~{a}o 1010,
Cidade Universit\'{a}ria, \\CEP 05508-090, S\~{a}o Paulo, S.P., Brazil \\
\\ Ronaldo Garcia\\Instituto de Matem\'{a}tica e
Estat\'{\i}stica,\\Universidade Federal de Goi\'as,\\CEP 74001-970, Caixa Postal
131,\\Goi\^ania, GO, Brazil}


\begin{thebibliography}{99}

\bibitem {bh}
 \noindent {\sc M.  Berry } and {\sc J.  Hannay,}
{\em Umbilic points on Gaussian Random
Surfaces, } Jour. Phys. A,  {\bf 10},  1977, 1809-1821.

\bibitem {bf}
 \noindent {\sc B. Bruce} and {\sc D. Fidal,} {\em On binary differential
equations and umbilic
points,} Proc. Royal Soc. Edinburgh  {\bf 111A}, 1989, pp. 147-168.


\bibitem {cj}
 \noindent {\sc  R. Courant} and {\sc F. John}, 
 {Introduction to Calculus and Analysis,
Vol. II},
  New York: Springer, 1989.



\bibitem {da}
 \noindent {\sc  G. Darboux, } {Le\c cons sur la Th\'eorie des Surfaces, vol. IV.
Sur la forme des lignes de courbure dans la voisinage d'un ombilic ,
Note 07},
  Paris:Gauthier Villars, 1896.

 \bibitem {ga}
\noindent{\sc  F.R.  Gantmacher, } Matrix Theory, Vol. I, Chelsea, 1960.


\bibitem{garcia} {\sc R. Garcia},   {\em Principal Curvature Lines near Darbouxian
 Partially Umbilic Points of Hypersurfaces Immersed in $\mathbb R^4$},
 Computational and Applied Mathematics, {\bf 20}, 2001, pp. 121-148.

\bibitem {a1}
   \noindent{\sc  R. Garcia   } and {\sc J.  Sotomayor,} {\em Structural stability of parabolic points and periodic asymptotic lines},
 Matem\'atica Contempor\^anea, {  \bf 12},   1997,  83-102.


\bibitem{ggs} {\sc  R. Garcia}, {\sc C. Gutierrez}   and {\sc J. Sotomayor}, { \em Structural
Stability of Asymptotic Lines on Surfaces Immersed in $\mathbb R^3$,} 
  Bulletin de Sciences Math\'ematiques, {\bf 123}, 1999, pp. 599-622.

\bibitem {m}
   \noindent{\sc R. Garcia } and {\sc  J. Sotomayor,} {\em Structurally stable configurations of lines of mean curvature and  umbilic points on surfaces immersed in ${\mathbb R}^3$,}    Publ. Matem\'atiques.  { \bf 45},  2001,   pp. 431-466.

\bibitem {g}
   \noindent{\sc R. Garcia } and {\sc  J. Sotomayor,} {\em Lines of Geometric Mean Curvature on surfaces immersed in ${\mathbb R}^3$,} to appear in 
    Annales de la Facult\'e des Sciences de Toulouse, {\bf 11}, 2002,  arXiv:math.DS/0302194 v1, 2003.
 

\bibitem {h}
   \noindent{\sc R. Garcia } and {\sc  J. Sotomayor,} {\em Lines of Harmonic Mean Curvature on surfaces immersed in ${\mathbb R}^3$,}
   To appear in Boletim da Soc. Bras. de Matem\'atica ( Bull. Bras. Math. Soc.), {\bf 34}, 2003, arXiv:math.DS/0302196 v1, 2003.
 

 

\bibitem{ax}\noindent{\sc R. Garcia } and {\sc  J. Sotomayor,} { \em Lines of Axial
Curvature on Surfaces Immersed into $\mathbb R^4$,}
 Differential Geometry and its Applications, {\bf 12}, 2000, 
pp. 253-269.
   
 
 \bibitem {gs1}
  \noindent{\sc  C. Gutierrez } and {\sc  J. Sotomayor,} 
 {\em Structurally Stable Configurations of  Lines of Principal  Curvature},
  Asterisque {\bf  98-99}, 1982, 185-215.

\bibitem{gsln} {\sc C. Gutierrez} and {\sc J. Sotomayor},
 {\em  An Approximation Theorem for Immersions with Structurally Stable
Configurations of
Lines of Principal Curvature,}
  Lect. Notes in Math. {\bf 1007},
  1983.

\bibitem {gs2}
  \noindent{\sc  C. Gutierrez  } and {\sc J. Sotomayor,} {Lines of
Curvature and Umbilic Points on Surfaces}, $18^{\hbox{th}}$  Brazilian
Math.
Colloquium, Rio de Janeiro,  IMPA, 1991.   Reprinted with update as  { Structurally Stable
Configurations of Lines of Curvature and Umbilic Points on Surfaces,
Lima, Monografias del IMCA}, 1998.

\bibitem{gggt} {\sc C. Gutierrez}, {\sc I. Guadalupe}, {\sc R. Tribuzy} and {\sc V. Gu\'\i \~nez},   
{\em Lines of curvature on surfaces immersed in $\mathbb R^4$},   Bol. Soc. Bras. Mat.    
{\bf 28}, 1997,  pp. 233-251.
 
  
\bibitem {iv}
 {\sc V. V. Ivanov},  {\em The  analytic conjecture of Carath\'eodory},  
 Siberian Math. Journal, {\bf 43}, 2002,   251--322.
 


\bibitem {ms}
 {\sc L.F. Mello } and {\sc J. Sotomayor},  {\em A note on some developments on Carath\'eodory Conjecture on umbilic points }, Expo. Math. {\bf 17},
  1999, 49-58.






\bibitem{p}
  {\sc I. R. Porteous, } Geometric  Differentiation, Cambridge Univ. Press,   1994.


\bibitem{sx} { \sc B. Smyth } and {\sc F. Xavier }
 {\em A sharp geometric estimate for the index of an umbilic point on a smooth surface,} 
   Bull. London Math. Soc., {\bf 24}, 1992, 176-180.

\bibitem{sp}
  {\sc M. Spivak, } Introduction to
Comprehensive Differential Geometry, Vol. III
   Berkeley, Publish or Perish, 1980.



\bibitem{st} {\sc D. Struik,}
  Lectures  on  Classical  Differential
 Geometry,
 Addison Wesley Pub. Co.,  Reprinted by Dover Publications,
Inc., 1988.

\bibitem{t} {\sc R. Thom,} Stabilit\'e Structurelle et Morphog\'en\`ese,
  W. A. Benjamin, Inc., 1972.


\bibitem {w}
  \noindent{\sc  T.C.T. Wall,}  
{\em Geometric properties of generic differentiable manifolds,} Springer Lecture Notes in Mathematics, {\bf 597}, 1977, 707-774.

\bibitem {z}
\noindent{\sc  E. C. Zeeman},  {\em The umbilic bracelet and the double cusp catastrophe}, Springer Lecture Notes in Mathematics, {\bf 525}, 328-366.



\end{thebibliography}
\end{document}